\documentclass[10pt, letterpaper]{article}
\usepackage{verbatim}
\usepackage{amsmath}
\usepackage{amssymb}
\usepackage{mathrsfs}
\usepackage{epsfig}
\usepackage{url}
\usepackage{hyperref}
\usepackage{graphicx}
\usepackage{float}
\usepackage{theorem}
\usepackage{color}

\usepackage{algorithm}
\usepackage{algorithmic}
\newcommand{\be}[1]{\begin{equation}\label{#1}}
\newcommand{\benon}{\begin{equation*}}  
\newcommand{\bemuln}[1]{\begin{multline}\label{#1}}
\newcommand{\bemul}{\begin{multline*}}
\newcommand{\bee}{\begin{eqnarray*}}
\newcommand{\eee}{\end{eqnarray*}}
\newcommand{\been}[1]{\begin{eqnarray}\label{#1}}
\newcommand{\eeen}{\end{eqnarray}}
\newcommand{\began}[1]{\begin{gather}\label{#1}}
\newcommand{\bega}{\begin{gather*}}
\newcommand{\bealn}[1]{\begin{align}\label{#1}}
\newcommand{\beal}{\begin{align*}}
\newcommand{\bealatn}[2]{\begin{alignat}{#1}\label{#2}}
\newcommand{\bealat}{\begin{alignat*}}
\newcommand{\bexalatn}[1]{\begin{xalignat}\label{#1}}
\newcommand{\bexalat}{\begin{xalignat*}}
\newcommand{\qed}{\newline \mbox{ } \hfill 
            \rule[-1pt]{2.5mm}{2.5mm}\par\vskip 10pt }

\newcommand{\pf}{\vskip 5 pt \noindent {\bf{Proof: }}}

\newcommand{\mb}{\mathbf}

\newcommand{\mbb}{\mathbb}

{\theoremstyle{plain} \newtheorem{thm}{Theorem}[section]

\newtheorem{lemma}[thm]{Lemma}

}
{\theoremstyle{break} \theorembodyfont{\it}
\newtheorem{defi}{Definition}
 
 }

\def\bv{{\mathbf v}}

\def\bx{{\mathbf x}}  
\def\by{{\mathbf y}}
\def\bz{{\mathbf z}}
\def\bA{{\mathbf A}}

\def\bE{{\mathbf E}}

\def\bI{{\mathbf I}}

\def\texitem#1{\par\smallskip\noindent\hangindent 25pt
               \hbox to 25pt {\hss #1 ~}\ignorespaces}

\newcommand{\bzero}{{\mathbf{0}}}

\newcommand{\scrA}{\mathcal{A}}

\newcommand{\scrD}{\mathcal{D}}

\newcommand{\scrF}{\mathcal{F}}

\newcommand{\scrS}{\mathcal{S}}

\newcommand{\scrX}{\mathcal{X}}

\newcommand{\btheta}{\boldsymbol{\theta}}

\newcommand{\bmu}{\boldsymbol{\mu}} 

\newcommand{\bpi}{{\boldsymbol{\pi}}}

\newcommand{\bphi}{{\boldsymbol{\phi}}}

\newcommand{\BS}{\boldsymbol}
\newcommand{\T}{\btheta}
\newcommand{\real}{{\mathbb R}}

\begin{document}
\title{Learning Policies for Markov Decision Processes from
  Data \thanks{Submitted for
    publication. Research partially supported by the NSF under grants
    CNS-1239021, CCF-1527292, and IIS-1237022, and by the ARO under
    grants W911NF-11-1-0227 and W911NF-12-1-0390. Hao Liu has been
    supported by the Lin Guangzhao \& the Hu Guozan Graduate Education
    International Exchange Fund.}}

\author{Manjesh~K.~Hanawal,\thanks{$\dagger$ M. K. Hanawal is with
    Industrial Engineering and Operations Research, IIT-Bombay, Powai,
    MH 400076, India. E-mail: {\tt mhanawal@iitb.ac.in.}}\ Hao~Liu,
  Henghui~Zhu,\thanks{$\ddagger$ H. Liu, H. Zhu are with the Center for
    Information and Systems Engineering, Boston University, Boston, MA
    02215, USA. H. Liu is also with the College of Control Science and
    Engineering, Zhejiang University, Hangzhou, Zhejiang, 310027
    China. E-mail: \texttt{\{haol, henghuiz\}@bu.edu}.}\\ and
  Ioannis~Ch.~Paschalidis \thanks{$\S$ Corresponding
    author. I. Ch. Paschalidis is with the Department of Electrical and
    Computer Engineering and the Division of Systems Engineering, Boston
    University, Boston, MA 02215, USA. E-mail: {\tt yannisp@bu.edu},
    URL: \url{http://sites.bu.edu/paschalidis/}.}  }

\markboth{IEEE Transactions on Automatic Control}%
{Hanawal \MakeLowercase{\textit{et al.}}: Learning Policies for Markov
  Decision Processes from Data}
\maketitle

\begin{abstract}
  We consider the problem of learning a policy for a Markov decision
  process consistent with data captured on the state-actions pairs
  followed by the policy. We assume that the policy belongs to a class
  of parameterized policies which are defined using features associated
  with the state-action pairs. The features are known a priori, however,
  only an unknown subset of them could be relevant. The policy
  parameters that correspond to an observed target policy are recovered
  using $\ell_1$-regularized logistic regression that best fits the
  observed state-action samples. We establish bounds on the difference
  between the average reward of the estimated and the original policy
  (regret) in terms of the generalization error and the ergodic
  coefficient of the underlying Markov chain. To that end, we combine
  sample complexity theory and sensitivity analysis of the stationary
  distribution of Markov chains. Our analysis suggests that to achieve
  regret within order $O(\sqrt{\epsilon})$, it suffices to use training
  sample size on the order of $\Omega(\log n \cdot poly(1/\epsilon))$,
  where $n$ is the number of the features. We demonstrate the
  effectiveness of our method on a synthetic robot navigation example.
\end{abstract}

%\begin{IEEEkeywords}
%Machine learning, Markov Decision Processes, reinforcement learning,
%regression.
%\end{IEEEkeywords}

% \fi
%\IEEEpeerreviewmaketitle

\section{Introduction}

Markov Decision Processes (MDPs) offer a framework for
many dynamic optimization problems under
uncertainty~\cite{puterman1994mdp,bert-dp}. When the state-action space
is not large and transition probabilities for all state-action pairs
are known, standard techniques such as policy iteration and value
iteration can compute an optimal policy. More often than not, however,
problem instances of realistic size have very large state-action spaces
and it becomes impractical to know the transition probabilities
everywhere and compute an optimal policy using these {\em off-line}
methods.

For such large problems, one resorts to approximate methods,
collectively referred to as {\em Approximate Dynamic Programming
  (ADP)}~\cite{beti96,pow-adp07}. ADP methods approximate either the
value function and/or the policy and optimize with respect to the
parameters, as for instance is done in {\em actor-critic}
methods~\cite{kots03,wapa-tac-hac-2016,pa-est-li-nrl-2012}.  The
optimization of approximation parameters requires the learner to have
access to the system and be able to observe the effect of applied
control actions.

In this paper, we adopt a different perspective and assume that the
learner has no direct access to the system but only has samples of the
actions applied in various states of the system. These actions applied
in various states are generated according to a policy that is fixed but
unknown. As an example, the actions could be followed by an expert
player who plays an optimal policy. Our goal is to learn a policy
consistent with the observed states-actions, which we will call {\em
  demonstrations}. 

Learning from an expert is a problem that has been studied in the
literature and referred to as apprenticeship
learning~\cite{ICML04_ApprenticeshipLearning_AbbeelNg,
  UAI2007_ApprenticeshipLearning_NeuSczepesvari}, imitation
learning~\cite{AISTATS2011_AReductionOfImitationLearning_RossGordonBagnell}
or learning from
demonstrations~\cite{ECML2010_LearningFromDemonstration_MeloLopes}. While
there are many settings where it could be useful, the main application
driver has been robotics~\cite{TCS99_ImitationLearning_Schaal,
  RAS09_ASurveyOfRobotLearning_ArgallChernovaVeloso, erlhagen2006goal}.
Additional interesting application examples include: learning from an
experienced human pilot/driver to navigate vehicles autonomously,
learning from animals to develop bio-inspired policies, and learning
from expert players of a game to train a computer player. In all these
examples, and given the size of the state-action space, we will not
observe the actions of the expert in all states, or more broadly
``scenarios'' corresponding to parts of the state space leading to
similar actions. Still, our goal is to learn a policy that generalizes
well beyond the scenarios that have been observed and is able to select
appropriate actions even at unobserved parts of the state space.

A plausible way to obtain such a policy is to learn a mapping of the
states-actions to a lower dimensional space. Borrowing from ADP methods,
we can obtain a lower-dimensional representation of the state-action
space through the use of features that are functions of the state and
the action taken at that state. In particular, we will consider policies
that combine various features through a weight vector and reduce the
problem of learning the policy to learning this weight/parameter
vector. Essentially, we will be learning a parsimonious parametrization
of the policy used by the expert.

The related work in the literature on learning through demonstrations
can be broadly classified into {\em direct} and {\em indirect}
methods~\cite{RAS09_ASurveyOfRobotLearning_ArgallChernovaVeloso}. In the
direct methods, supervised learning techniques are used to obtain a best
estimate of the expert's behavior; specifically, a best fit to the
expert's behavior is obtained by minimizing an appropriately defined
loss function. A key limitation of this method is that the estimated
policy is not well adapted to the parts of the state space not visited
often by the expert, thus resulting in poor action selection if the
system enters these states. Furthermore, the loss function can be
non-convex in the policy, rendering the corresponding problem hard to
solve.

Indirect methods, on the other hand, evaluate a policy by learning the
full description of the MDP. In particular, one solves a so called {\em
  inverse reinforcement learning} problem which assumes that the
dynamics of the environment are known but the one-step reward function
is unknown. Then, one estimates the reward function that the expert is
aiming to maximize through the demonstrations, and obtains the policy
simply by solving the MDP. A policy obtained in this fashion tends to
generalize better to the states visited infrequently by the expert,
compared to policies obtained by direct methods. The main drawback of
inverse reinforcement learning is that at each iteration it requires to
solve an MDP which is computationally expensive. In addition, the
assumption that the dynamics of the environment are known for all states
and actions is unrealistic for problems with very large state-action
spaces.

In this work, we exploit the benefits of both direct and indirect
methods by assuming that the expert is using a {\em Randomized
  Stationary Policy (RSP)}. As we alluded to earlier, an RSP is
characterized in terms of a vector of features associated with
state-action pairs and a parameter $\T$ weighing the various elements of
the feature vector. We consider the case where we have many features,
but only relatively few of them are sufficient to approximate the target
policy well. However, we do not know in advance which features are
relevant; learning them and the associated weights (elements of $\T$) is
part of the learning problem.

We will use supervised learning to obtain the best estimate of the
expert's policy. As in actor-critic methods, we use an RSP which is a
parameterized ``Boltzmann'' policy and rely on an $\ell_1$-regularized
maximum likelihood estimator of the policy parameter vector. An
$\ell_1$-norm regularization induces sparse estimates and this is useful
in obtaining an RSP which uses only the relevant features. In
\cite{ICML2004_FeatureSelection_Ng}, it is shown that the sample
complexity of $\ell_1$-penalized logistic regression grows as
$\mathcal{O}(\log n)$, where $n$ is the number of features. As a result,
we can learn the parameter vector $\T$ of the target RSP with relatively
few samples, and the RSP we learn generalizes well across states that
are not included in the demonstrations of the expert. Furthermore,
$\ell_1$-regularized logistic regression is a convex optimization
problem which can be solved efficiently.

\subsection{Related Work}
There is substantial work in the literature on learning MDP policies by
observing experts; see
\cite{RAS09_ASurveyOfRobotLearning_ArgallChernovaVeloso} for a
survey. We next discuss papers that are more closely related to our
work.

In the context of indirect methods,
\cite{ICML2000_AlgorithmsForInverse_NgRussel} develops techniques for
estimating a reward function from demonstrations under varying degrees
of generality on the availability of an explicit
model. \cite{ICML04_ApprenticeshipLearning_AbbeelNg} introduces an
inverse reinforcement learning algorithm that obtains a policy from
observed MDP trajectories followed by an expert. The policy is
guaranteed to perform as well as that of the expert's policy, even
though the algorithm does not necessarily recover the expert's reward
function. In \cite{UAI2007_ApprenticeshipLearning_NeuSczepesvari}, the
authors combine supervised learning and inverse reinforcement learning
by fitting a policy to empirical estimates of the expert demonstrations
over a space of policies that are optimal for a set of parameterized
reward functions. \cite{UAI2007_ApprenticeshipLearning_NeuSczepesvari}
shows that the policy obtained generalizes well over the entire state
space. \cite{ECML2010_LearningFromDemonstration_MeloLopes} uses a
supervised learning method to learn an optimal policy by leveraging the
structure of the MDP, utilizing a kernel-based approach. Finally,
\cite{AISTATS2011_AReductionOfImitationLearning_RossGordonBagnell}
develops the DAGGER (Dataset Aggregation) algorithm that trains a
deterministic policy which achieves good performance guarantees under
its induced distribution of states.

Our work is different in that we focus on a parameterized set of
policies rather than parameterized rewards. In this regard, our work is
similar to approximate DP methods which parameterize the policy, e.g.,
expressing the policy as a linear functional in a lower dimensional
parameter space. This lower-dimensional representation is critical in
overcoming the well known ``curse of dimensionality.''

\subsection{Contributions}
We adopt the $\ell_1$-regularized logistic regression to estimate a
target RSP that generates a given collection of state-action
samples. We evaluate the performance of the estimated policy and derive
a bound on the difference between the average reward of the estimated
RSP and the target RSP, typically referred to as {\em regret}.  We show
that a good estimation of the parameter of the target RSP also implies a
good bound on the regret. To that end, we generalize a sample complexity
result on the log-loss of the maximum likelihood estimates
\cite{ICML2004_FeatureSelection_Ng} from the case of two actions
available at each state to the multi-action case. Using this result, we
establish a sample complexity result on the regret.

Our analysis is based on the novel idea of separating the loss in
average reward into two parts. The first part is due to the error in the
policy estimation (training error) and the second part is due to the
perturbation in the stationary distribution of the Markov chain caused
by using the estimated policy instead of the true target policy
(perturbation error). We bound the first part by relating the result on
the log-loss error to the Kulback-Leibler divergence
\cite{Book_ElementsOfInformation} between the estimated and the target
RSPs. The second part is bounded using the ergodic coefficient of the
induced Markov chain. Finally, we evaluate the performance of our method
on a synthetic example.

The paper is organized as follows. In Section \ref{sec: Setup}, we
introduce some of our notation and state the problem. In Section
\ref{sec: Regression}, we describe the supervised learning algorithm
used to train the policy. In Section~\ref{sec:log-loss}, we establish a
result on the (log-loss) error in policy estimation. In
Section~\ref{sec:regret}, we establish our main result which is a bound
on the regret of the estimated policy. In Section \ref{sec: Simulation},
we introduce a robot navigation example and present our numerical
results. We end with concluding remarks in Section \ref{sec:Conclusion}.

\textbf{Notational conventions}. Bold letters are used to denote vectors
and matrices; typically vectors are lower case and matrices upper
case. Vectors are column vectors, unless explicitly stated
otherwise. Prime denotes transpose.  For the column vector $\bx\in
\mbb{R}^n$ we write $\bx=(x_1,\ldots,x_n)$ for economy of space, while
$\|\bx\|_p=(\sum_{i=1}^{n}|x_i|^p)^{1/p}$ denotes its $p$-norm.  Vectors
or matrices with all zeroes are written as $\bzero$, the identity matrix
as $\bI$, and $\mb{e}$ is the vector with all entries set to $1$.  For
any set $\scrS$, $|\scrS|$ denotes its cardinality. We use $\log$ to
denote the natural logarithm and a subscript to denote different bases,
e.g., $ \log_2$ denotes logarithm with base $2$.

\section{Problem Formulation}
\label{sec: Setup}

Let $(\mathcal{X}, \mathcal{A}, \BS{P}, R)$ denote a Markov Decision
Process (MDP) with a finite set of states $\mathcal{X}$ and a finite set
of actions $\mathcal{A}$.  For a state-action pair $(x, a) \in
\mathcal{X}\times \mathcal{A}$, let $P(y |x,a)$ denote the probability
of transitioning to state $y \in \mathcal{X}$ after taking action $a$ in
state $x$. The function $R: \mathcal{X}\times \mathcal{A} \rightarrow
\mbb{R}$ denotes the one-step reward function. 

Let us now define a class of {\em Randomized Stationary Policies (RSPs)}
parameterized by vectors $\T \in \mbb{R}^n$. Let $\{\mu_{\T}: \T \in
\mbb{R}^n\}$ denote the set of RSPs we are considering. For a given
parameter $\T$ and a state $x \in \mathcal{X}$, $\mu_{\T}(a|x)$ denotes
the probability of taking action $a$ at state $x$. Specifically, we 
consider the Boltzmann-type RSPs of the form
\begin{equation}
\label{eqn:RSP}
\mu_{\T}(a|x)=\frac{\exp\{\T' \bphi(x,a)\}}{\sum_{b \in \mathcal{A}}
  \exp\{\T' \bphi(x,b)\}}, 
\end{equation}
where $\bphi: \mathcal{X}\times \mathcal{A} \rightarrow [0,\,1]^n$ is a
vector of features associated with each state-action pair
$(x,a)$. (Features are normalized to take values in $[0,1]$.)
Henceforth, we identify an RSP by its parameter $\T$. We assume that the
policy is {\em sparse}, that is, the vector $\T$ has only $r<n$ non-zero
components and each is bounded by $K$, i.e., $|\theta_i| < K$ for all
$i$. Given an RSP $\T$, the resulting transition probability matrix on
the Markov chain is denoted by $\BS{P}_{\T}$, whose $(x,y)$ element is
$P_{\T}(y|x)=\sum_{a\in \mathcal{A}}\mu_{\T}(a|x) P(y|x,a)$ for all
state pairs $(x,y)$.

Notice that for any RSP $\T$, the sequence of states $\{X_k\}$ and
the sequence of state-action pairs $\{X_k, A_k\}$ form a Markov chain with
state space $\mathcal{X}$ and $\mathcal{X}\times\mathcal{A}$,
respectively. We assume that for every $\T$, the Markov chains $\{X_k\}$
and $\{X_k, A_k\}$ are irreducible and aperiodic with stationary
probabilities $\pi_{\BS{\theta}}(x)$ and
$\eta_{\BS{\theta}}(x,a)=\pi_{\BS{\theta}}(x)\mu_{\BS{\theta}}(a|x)$,
respectively.

The average reward function associated with an RSP $\T$ is a function
$\overline{R}: \mbb{R}^n \rightarrow \mbb{R}$ defined as
\begin{equation}
\label{eqn:AverageReward}
\overline{R}(\T)=\sum_{(a,x)}\eta_{\T}(x,a)R(x,a).
\end{equation}

Let now fix a target RSP $\T^*$. As we assumed above, $\T^*$ is sparse
having at most $r$ non-zero components $\theta^*_i$, each satisfying
$|\theta^*_i|\leq K$. This is simply the policy used by an expert (not
necessarily optimal) which we wish to learn. We let
$\scrS(\T^*)=\{(x_i,a_i): i=1,2,\ldots,m\}$ denote a set of state-action
samples generated by playing policy $\T^*$. The state samples $\{x_i
:i=1,2,\ldots,m \}$ are independent and identically distributed (i.i.d.)
drawn from the stationary distribution $\bpi_{\T^*}$ and $a_i$ is the
action taken in state $x_i$ according to the policy $\mu_{\T^*}$. It
follows that the samples in $\scrS(\T^*)$ are i.i.d. according to the
distribution $\mathcal{D} \sim \eta_{\T^*}(x,a)$.

We assume we have access only to the demonstrations $\scrS(\T^*)$ while
the transition probability matrix $\BS{P}_{\T^*}$ and the target RSP
$\T^*$ are unknown. The goal is to learn the target policy $\T^*$ and
characterize the average reward obtained when the learned RSP is applied
to the Markov process. In particular, we are interested in estimating
the target parameter $\T^*$ from the samples efficiently and evaluate
the performance of the estimated RSP with respect to the target RSP,
i.e., bound the regret defined as
\[
\mathrm{Reg}(\scrS(\T^*))=
\overline{R}(\T^*)-\overline{R}(\hat{\btheta}), 
\] 
where $\hat{\btheta}$ is the estimated RSP from the samples $\scrS(\T^*)$.

\section{Estimating the policy} 
\label{sec: Regression}

Next we discuss how to estimate the target RSP $\T^*$ from the $m$
i.i.d. state-action training samples in $\scrS(\T^*)$. Given the
Boltzmann structure of the RSP we have assumed, we fit a logistic
regression function using a regularized maximum likelihood estimator as
follows:
\begin{equation}
\label{eqn:logistic}
\begin{array}{rl} 
\max_{\T \in \real^n}  & \sum_{i=1}^m \log \mu_{\T}(a_i|x_i)\\
\text{s.t.} & \|\T \|_1 \leq B,
\end{array}
\end{equation}
where $B$ is a parameter that adjusts the trade-off between fitting the
training data ``well'' and obtaining a sparse RSP that generalizes well
on a test sample.

We can evaluate how well the maximum likelihood function fits the
samples in the logistic function using a log-loss metric, defined as the
expected negative of the likelihood over the (random) test
data. Formally, for any parameter $\T$, log-loss is given by
\begin{equation} \label{log-loss}
\epsilon(\T)=\mathbb{E}_{(x,a)\sim \mathcal{D}}[-\log \mu_{\T}(a|x)],
\end{equation}
where the expectation is taken over state-action pairs $(x,y)$ drawn
from the distribution $\mathcal{D}$; recall that we defined
$\mathcal{D}$ to be the stationary distribution $\eta_{\T^*}(x,a)$ of
the state-action pairs induced by the policy $\T^*$. Since the
expectation is taken with respect to new data not included in the
training set $\scrS(\T^*)$, we can think of $\epsilon(\T)$ as an {\em
  out-of-sample} metric of how well the RSP $\T$ approximates the
actions taken by the target RSP $\T^*$.

We also define a sample-average version of log-loss: given any set
$\scrS=\{(x_i,a_i): i=1,2,\ldots,m\}$ of state-action pairs, define
\begin{equation} \label{sample-log-loss}
\hat{\epsilon}_{\scrS}(\T)=\frac{1}{m}  \sum_{i=1}^m (-\log
\mu_{\T}(a_i|x_i)).
\end{equation}
We will use the term {\em empirical} log-loss to refer to log-loss over
the training set $\scrS(\T^*)$ and use the notation
\be{empirical-log-loss} 
\hat{\epsilon}(\T)\stackrel{\triangle}{=}
\hat{\epsilon}_{\scrS(\T^*)}(\T).
\end{equation}

To estimate an RSP, say $\hat{\T}$, from the training set, we adopt the
logistic regression with $\ell_1$ regularization introduced in
\cite{ICML2004_FeatureSelection_Ng}.  The specific steps are shown in
Algorithm~\ref{alg:train}.
\begin{algorithm}[H]
\caption{Training algorithm to estimate the target RSP $\T^*$ from the
samples $\scrS(\T^*)$.}
\label{alg:train}
\begin{algorithmic}
\STATE {\em Initialization:} Fix $0<\gamma <1$ and $C \geq rK.$

\STATE Split the training set $\scrS(\T^*)$ into set two sets $\scrS_1$ and
$\scrS_2$ of size $(1-\gamma) m$ and $\gamma m$ respectively. $\scrS_1$
is used for training and $\scrS_2$ for cross-validation.

\STATE {\em Training:}
\FOR{$B=0,1,2,4,\ldots,C$} \STATE
Solve the optimization problem (\ref{eqn:logistic})  for each
$B$ on the set $\scrS_1$, and let $\T_B$ denote the optimal
solution.
\ENDFOR

\STATE {\em Validation:} Among the $\T_B$'s from the training step,
select the one with the lowest ``hold-out'' error on $\scrS_2$, i.e,
$\hat{B}=\arg \min_{B\in \{0,1,2,4\ldots,C\}}
\hat{\epsilon}_{\scrS_2}(\T_B)$ and set $\hat{\T}=\T_{\hat{B}}$.
\end{algorithmic}
\end{algorithm}

\section{Log-loss performance} \label{sec:log-loss}

In this section we establish a sample complexity result indicating that
relatively few training samples (logarithmic in the dimension $n$ of the
RSP parameter $\T$) are needed to guarantee that the estimated RSP
$\hat{\T}$ has out-of-sample log-loss close to the target RSP $\T^*$.
We will use the line of analysis in \cite{ICML2004_FeatureSelection_Ng}
but generalize the result from the case where only two actions are
available for every state to the general multi-action setting.

We start by relating the difference of the log-loss function associated
with RSP $\T$ and its estimate $\hat{\T}$ to the relative entropy, or
Kullback-Leibler (KL) divergence, between the corresponding RSPs. For a
given $x \in \mathcal{X}$, we denote the KL-divergence between RSPs,
$\T_1$ and $\T_2$ as $D(\BS{\mu}_{\T_1}(\cdot | x)\|
\BS{\mu}_{\T_2}(\cdot | x))$, where $\BS{\mu}_{\T}(\cdot|x)$ denotes the
probability distribution induced by RSP $\T$ on $\mathcal{A}$ in state
$x$, and is given as follows:
\[
D(\BS{\mu}_{\T_1}(\cdot | x)\|
\BS{\mu}_{\T_2}(\cdot | x))= \sum_a \mu_{\T_1}(a|x)\log \frac
{\mu_{\T_1}(a|x)}{\mu_{\T_2}(a|x)}.
\]

We also define the average KL-divergence, denoted by
$D_{\T}(\BS{\mu}_{\T_1}\| \BS{\mu}_{\T_2})$, as the average of
$D(\BS{\mu}_{\T_1}(\cdot | x)\| \BS{\mu}_{\T_2}(\cdot | x))$ over states
visited according to the stationary distribution $\BS{\pi}_{\T}$ of the
Markov chain induced by policy $\T$. Specifically, we define
\[
D_{\T}(\BS{\mu}_{\T_1}\| \BS{\mu}_{\T_2})=\sum_x
\pi_{\T}(x)D(\BS{\mu}_{\T_1}(\cdot|x)\| \BS{\mu}_{\T_2}(\cdot|x)).
\]

\begin{lemma}
Let $\hat{\T}$ be an estimate of RSP $\T$. Then,
\[
\epsilon(\hat{\T})-\epsilon(\T) = D_{\T}(\BS{\mu}_{\T}\|
\BS{\mu}_{\hat{\T}}) 
\] 
\end{lemma}
\pf 
\begin{eqnarray*}
  \lefteqn{\epsilon(\hat{\T})-\epsilon(\T)}\\
  &=& \sum_{x \in \mathcal{X},a \in \mathcal{A}}  \eta_{\T}(x,a)[\log
  \mu_{\T}(a|x) -\log \mu_{\hat{\T}}(a|x)]\\ 
  &=& \sum_{x\in \scrX} \pi_{\T}(x)\sum_{a\in \scrA} \mu_{\T}(a|x)\log
  \frac {\mu_{\T}(a|x)}{\mu_{\hat{\T}}(a|x)}\\ 
  &=& \sum_{x\in \scrX} \pi_{\T}(x)D(\BS{\mu}_{\T}(\cdot|x)\|
  \BS{\mu}_{\hat{\T}}(\cdot|x))\\
  &=& D_{\T}(\BS{\mu}_{\T}\| \BS{\mu}_{\hat{\T}}).
\end{eqnarray*}
\qed

For the case of a binary action at every state,
\cite{ICML2004_FeatureSelection_Ng} showed the following result. To
state the theorem, let $poly(\cdot)$ denote a function which is
polynomial in its arguments and recall that $m$ is the number of
state-action pairs in the training set $\scrS(\T^*)$ used to learn
$\hat{\T}$.

\begin{thm}[\cite{ICML2004_FeatureSelection_Ng}, Thm. 3.1]
\label{thm:SampleComlexity}
Suppose action set $\mathcal{A}$ contains only two actions,
i.e. $\mathcal{A} = \{0, 1\}$. Let $\epsilon >0$ and $\delta >0$. In
order to guarantee that, with probability at least $1-\delta$,
$\hat{\T}$ produced by the algorithm performs as well as $\T^*$, i.e.,
\[
D_{\T^*}(\bmu_{\T^*}\| \bmu_{\hat{\T}})\leq \epsilon,
\]
it suffices that $m = \Omega((\log n) \cdot poly (r, K, C,
\log(1/\delta), 1/\epsilon))$. 
\end{thm}

We will generalize the result to the case when more than two actions are
available at each state. We assume that $|\mathcal{A}|=H$, that is, at
most $H$ actions are available at each state. By introducing in
(\ref{eqn:RSP}) features that get activated at specific states, it is
possible to accommodate MDPs where some of the actions are not available
at these states. 

\begin{thm}\label{thm:1}
  Let $\varepsilon > 0$ and $\delta > 0$. In order to guarantee that,
  with probability at least $1-\delta$, $\hat{\T}$ produced by
  Algorithm~\ref{alg:train} performs as well as $\T^*$, i.e.,
\begin{equation}
|\epsilon(\hat{\T}) - \epsilon(\T^*) |< \varepsilon,
\end{equation}
it suffices that
\begin{equation*}
m = \Omega\Big( (\log n) \cdot \textrm{poly} (r,\,K,\,C,\,\,H,\,\log(1/\delta),\,1/\varepsilon) \Big).
\end{equation*}
Furthermore, in terms of only $H$, $m=\Omega(H^3)$.
\end{thm}

The remainder of this section is devoted to proving Theorem~\ref{thm:1}.
The proof is similar to the proof of Theorem~\ref{thm:SampleComlexity}
(Theorem 3.1 in \cite{ICML2004_FeatureSelection_Ng}) but with key
differences to accommodate the multiple actions per state. We start by
introducing some notations and stating necessary lemmata.

Denote by $\mathcal{F}$ a class of functions over some domain
$\scrD_{\scrF}$ and range $[-M, M] \subset \real$.  Let
$\mathcal{F}|_{\bz^{(1)},\ldots,\bz^{(m)}} = \{
(f(\bz^{(1)}),\ldots,f(\bz^{(m)})) \mid f \in \mathcal{F}\}\subset [-M,
M]^m$, which is the valuation of the class of functions at a certain
collection of points $\bz^{(1)}, \ldots,\linebreak[3] \bz^{(m)} \in
\scrD_{\scrF}$. It 
is said that a set of vectors 
$\{\bv^{(1)},\ldots,\bv^{(k)}\}$ in $\real^m$  $\varepsilon$-covers
$\mathcal{F}|_{\bz^{(1)},\ldots,\bz^{(m)}}$ in the $p$-norm, if for
every point $\bv \in \mathcal{F}|_{\bz^{(1)},\ldots,\bz^{(m)}}$ there
exists some $\bv^{(i)}$, $i=1,\ldots,k$, such that $\|\bv -
\bv^{(i)}\|_p \leq m^{1/p} \varepsilon$. Let also denote
$\mathcal{N}_p(\mathcal{F},\varepsilon,(\bz^{(1)},\,\ldots,\,\bz^{(m)}))$
the size of the smallest set that $\varepsilon$-covers
$\mathcal{F}|_{\bz^{(1)},\ldots,\bz^{(m)}}$ in the $p$-norm. Finally,
let 
\[ 
\mathcal{N}_p(\mathcal{F},\varepsilon,m) =
\sup_{\bz^{(1)},\ldots,\bz^{(m)}}
\mathcal{N}_p(\mathcal{F},\varepsilon,(\bz^{(1)},\ldots,\bz^{(m)})).
\]

To simplify the representation of log-loss of the general logistic
function in \eqref{eqn:RSP}, we use the following notations. For each $x
\in \mathcal{X}$, $\bphi(x)= (\bphi_1(x),\ldots,\bphi_H(x))
\linebreak[3] =
(\bphi(x,1),\ldots,\bphi(x,H)) \in [0, 1]^{nH}$ denotes feature vectors
associated with each action. For any $\T \in \real^n$, define the
$\log$-likelihood function $l: [0, 1]^{nH} \times \{1,\ldots,H\} \to
\real$ as
\begin{equation*} 
l(\bphi(x),i) =
-\log \left(
  \frac{\exp(\T'\bphi_i(x))}{\sum_{j=1}^H\exp(\T'\bphi_j(x))} \right). 
\end{equation*}
Note $\mu_{\T}(a|x)=\exp\{-l(\bphi(x),a)\}$. Further, let $g:[0, 1]^n
\rightarrow \real$ be the class of functions $g(\bx)=\T'\bx$. We can
then rewrite $l$ using $g$ as $l(\bphi(x),y)=l(g(\bphi_1(x)),
\linebreak[3] \ldots, 
g(\bphi_H(x)),y)$.

\begin{lemma}[\cite{ICML2004_FeatureSelection_Ng,pollard2012convergence}]
\label{lem:1} 
Let there be some distribution $D$ over $\scrD_{\scrF}$ and suppose
$\bz^{(1)},\ldots,\bz^{(m)}$ are drawn from $D$ i.i.d. Then,
\begin{multline*}
  \mbb{P} \left[ \exists f \in \mathcal{F} : \left| \frac{1}{m}
      \sum_{i=1}^{m} f(\bz^{(i)}) - \mbb{E}_{\bz \sim D} [f(\bz)]
    \right| \geq
    \varepsilon \right] \\
  \leq 8 \mbb{E}\big[\mathcal{N}_1(\mathcal{F}, \varepsilon/8,
  (\bz^{(1)},\ldots,\bz^{(m)})) \big] \exp\left(
  \frac{-m\varepsilon^2}{512 M^2}\right). 
\end{multline*}
\end{lemma}

\begin{lemma}[\cite{Zhang2002}]\label{lem:2}
  Suppose $\mathcal{G} = \{g: g(\bx) = \T' \bx,\ \bx \in \real^n,\
  \|\T\|_q \leq B \}$ and the input $\bx \in \real^n$ has a norm-bound
  such that $\|\bx\|_p \leq \zeta$, where $1/p+1/q = 1$. Then
\begin{equation}
\log_2 \mathcal{N}_2(\mathcal{G},\varepsilon,m) \leq
	\frac{B^2 \zeta^2}{\varepsilon^2}  \log_2(2n+1).
 \label{eq:lem2}
\end{equation}
\end{lemma}

\begin{lemma}[\cite{ICML2004_FeatureSelection_Ng,anthony2009neural}]
\label{lem:4} 
If $|f(\btheta) - \hat{f}(\btheta)| \leq \varepsilon$ for all $\btheta
\in \Theta$, then
\begin{equation*}
f\Big(\arg\min_{\btheta \in \Theta} \hat{f}(\btheta) \Big) \leq 
\min_{\btheta \in \Theta} f(\btheta) + 2\varepsilon.
\end{equation*}
\end{lemma}

\begin{lemma}\label{lem:modif}
  Let $\mathcal{G}$ be a class of functions from $\real^n$ to some
  bounded set in $\real$. Consider $\mathcal{F}$ as a class of functions
  from $\real^H \times \scrA$ to some bounded set in $\real$, with the
  following form:
\begin{multline}
\mathcal{F} = \big\{
f_g(\bphi(x),y) = l(g(\bphi_1(x)),\ldots,
g(\bphi_H(x)), y), g \in \mathcal{G},\ y \in\mathcal{A} \big\}.
\label{classF}
\end{multline}
If $l(\cdot,y)$ is Lipschitz with Lipschitz constant $L$ in the
$\ell_1$-norm for every $y \in \scrA$, then we have
\begin{equation*}
  \mathcal{N}_1(\mathcal{F},\varepsilon,m) \leq
  \big[
  \mathcal{N}_1(\mathcal{G}, \varepsilon /(LH), m)
  \big]^H.
\end{equation*}
\end{lemma}
\pf
Let $\Gamma = \mathcal{N}_1(\mathcal{G}, \varepsilon/(LH), m)$. It is
sufficient to show for every $m$ inputs 
\[ 
\bz^{(1)}= (\bphi(x^{(1)}),y^{(1)}), \ldots,
\bz^{(m)}=(\bphi(x^{(m)}),y^{(m)})
\]
we can find $\Gamma^H$ points in $\real^m$ that $\varepsilon$-cover
$\mathcal{F}|_{\bz^{(1)},\ldots,\bz^{(m)}}$.  Fix some set of points
$\bz^{(1)},\ldots,\bz^{(m)} \subset \real^{nH+1}$. From the definition
of $\mathcal{N}_1(\mathcal{G}, \varepsilon/(LH), m)$, for each $j \in
\{1,\ldots,H\}$, 
\begin{equation*}
\mathcal{N}_1(\mathcal{G}, \varepsilon/(LH),
(\bphi_j(x^{(1)}),\ldots,\bphi_j(x^{(m)}))) \leq \Gamma. 
\end{equation*}
Let~\footnote{One may find less than $\Gamma$ points, but we consider
  the worst case scenario.} $\{\bv_j^{(1)},\ldots,\bv_j^{(\Gamma)} \}$
be a set of $\Gamma$ points in $\real^m$
that $\varepsilon/(LH)$-covers \linebreak[4]
$\mathcal{G}|_{\bphi_j(x^{(1)}),\ldots,\bphi_j(x^{(m)})}$.  We use
notation $v_{j,i}^{k}$ to denote the $i$th element of vector
$\bv_{j}^{k}$. Then, for any $g \in \mathcal{G}$ and $j \in \{1,\ldots,
H\}$, there exists a $k(j) \in \{1,\ldots,\Gamma\}$ such that
\begin{align}
\label{eqn:CoverOfG}
\| (g(\bphi_j(x^{(1)})), & \ldots,g(\bphi_j(x^{(m)})) -
\bv_{j}^{k(j)}\|_1 \notag \\
= & \sum_{i=1}^m |g(\bphi_j(x^{(i)})) - v_{j,i}^{k(j)}| \notag \\
\leq & m \frac{\varepsilon}{LH}.
\end{align}

Now consider $\Gamma^H$ points with the following form
\bemul
l( v_{1,1}^{(j_1)},v_{2,1}^{(j_2)},\ldots,v_{H,1}^{(j_H)},y^{(1)}),
\ldots,%\\ 
l( v_{1,m}^{(j_1)},v_{2,m}^{(j_2)},\ldots,v_{H,m}^{(j_H)},y^{(m)} ),
\end{multline*}
where $j_1,\ldots,j_H \in \{1,\ldots, \Gamma\}$.

Given a $g \in \mathcal{G}$ and $f_g(\cdot) \in \mathcal{F}$, let $k(1),
\ldots k(H) \in \{1,\ldots, \Gamma\}$ be as defined above and consider
\begin{align*}
& \bigg\| \big(f_g(\bz^{(1)}),\ldots, f_g(\bz^{(m)})\big) - \\
& \hspace{1cm} \big( l(v_{1,1}^{k(1)},v_{2,1}^{k(2)},\ldots,
  v_{H,1}^{k(H)},y^{(1)}), \ldots, \\
& \hspace{1cm}  l(v_{1,m}^{k(1)}, v_{2,m}^{k(2)}, \ldots,
v_{H,m}^{k(H)}, y^{(m)}) \big) \bigg\|_1 \\
& \leq 
L \bigg\| \bigg( \sum_{h=1}^{H} |v_{h,1}^{k(h)} -g(\bphi_h(x^{(1)}))|,
\ldots,\\ 
& \hspace{1cm} \sum_{h=1}^{H} |v_{h,m}^{k(h)}
-g(\bphi_h(x^{(m)}))|\bigg) \bigg\|_1\\ 
& = L \sum_{h=1}^{H} \sum_{i=1}^m |g(\bphi_h(x^{(i)})) -
v_{h,i}^{k(h)}| \leq  m \varepsilon, 
\end{align*}
where we used the Lipschitz property of the $l(\cdot)$ function in the
first inequality and the last inequality follows from
(\ref{eqn:CoverOfG}). Thus, the $\Gamma^H$ points $\epsilon$-cover
$\mathcal{F}|_{\bz^{(1)},\ldots,\bz^{(m)}}$ in $\ell_1$-norm. Finally,
notice that the set of $m$ points $\{\bz^{(1)},\ldots,\bz^{(m)}\}$ is
arbitrary, which concludes the proof. \qed

We now continue with the proof of Theorem \ref{thm:1}. Recall
Algorithm~\ref{alg:train} and let $\hat{B}$ be the smallest integer in
$\{0,1,2,4,\ldots\}$ that is greater or equal to $rK$. Notice that in
Algorithm~\ref{alg:train} one can use a larger $C$ but we select the
smallest possible to obtain a tighter bound. For such a $\hat{B}$, it
follows that $rK\leq \hat{B} \leq 2rK$. Define a class of functions
$\mathcal{G}$ with domain $[0, 1]^n$ as
\begin{equation*}
\mathcal{G} = \Big\{ g: [0,1]^n \to \real \; \Big|\;  g(\bx) = \T'
\bx,\, \| \T \|_1 \leq \hat{B} \Big\}. 
\end{equation*}
By Lemma \ref{lem:2} and Eq. \eqref{eq:lem2},
\begin{equation*}
  \log_2 \mathcal{N}_2(\mathcal{G},\varepsilon/H,m) \leq 
  \frac{\hat{B}^2 H^2}{\varepsilon^2} \log_2(2n+1). 
\end{equation*}
The partial derivatives of the log-loss function are
\begin{equation*}
\frac{\partial}{\partial x_i} l(x_1,\ldots,x_H,k) =
\begin{cases}
- 1 + \frac{\exp(x_i)}{\sum_{j = 1}^{H}\exp(x_j) },\; k = i,\\
\frac{\exp(x_i)}{\sum_{j = 1}^{H}\exp(x_j) },\; k \neq i,
\end{cases}
\end{equation*}
and it can be seen that 
\[ 
\left|\frac{\partial}{\partial x_i}
l(x_1,\ldots,x_H,k)\right| \leq 1.
\]  
Hence, the Lipschitz constant for $l(\cdot,y)$ is $1$ for any $y =
1,\ldots,H$.

By Lemma \ref{lem:modif}, we have
\begin{equation*}
\log_2 \mathcal{N}_1(\mathcal{F}, \varepsilon, m)
\leq
H \log_2 \mathcal{N}_1 (\mathcal{G}, \varepsilon/H, m).
\end{equation*}
Using the relation $\mathcal{N}_1\leq \mathcal{N}_2$
(\cite{ICML2004_FeatureSelection_Ng, Zhang2002, anthony2009neural}), we
obtain
\begin{equation}
\log_2 \mathcal{N}_1(\mathcal{F}, \varepsilon, m)
\leq
\frac{\hat{B}^2 H^3}{\varepsilon^2} \log_2(2n+1). \label{eq:n1}
\end{equation}

We next find the range of class $\mathcal{F}$. To begin with, the range
of class $\mathcal{G}$ is
\begin{equation*}
| g(\bx) | = | \T' \bx| \leq \|\T\|_1 \|\bx\|_\infty \leq \hat{B}.
\end{equation*}
Since $l(\cdot,i)$ is Lipschitz in $\ell_1$-norm with Lipschitz constant
$1$ and $|f(0,\ldots,0,y)|\linebreak[3] =\log(H) < H$ (by the fact $H \geq 2$), then
\[
|f_g(\bphi(x),y)-f(\bzero,y)| \leq \sum_{h=1}^H |\T'\bphi_h(x)| \leq
H\hat{B},
\]  
which implies
\begin{equation}
|f_{g}(\bphi(x),y)| \leq H\hat{B} + H. \label{eq:m}
\end{equation}

Finally, let $m_1 = (1-\gamma) m$, which is the size of training set in
Algorithm~\ref{alg:train}. From Lemma~\ref{lem:1}, Eq.~\eqref{eq:m} and
Eq.~\eqref{eq:n1}, we have
\begin{multline}
 \mbb{P} \bigg[ \exists f \in \mathcal{F}:\; \bigg| \frac{1}{m_1}
\sum_{i=1}^{m_1} f(\bphi(x^{(i)}),y^{(i)}) - \mbb{E}_{\bz\sim D} [f(\bz)]
\bigg| \geq \varepsilon \bigg]\\
\leq  8 \cdot 2^{\frac{256 r^2 K^2 H^3}{\varepsilon^2}} (2n+1) \exp\bigg( \frac{-m_1 \varepsilon^2}{512(2rK + 1)^2H^2} \Bigg).
\end{multline} 

Treat $(1-\gamma)$ as a constant. To upper bound the right hand side of
the above equation by $\delta$, it suffices to have
\begin{equation}
m = \Omega\Big( (\log n) \cdot \textrm{poly} (r,\,K,\,H,\,\,\log(1/\delta),\,1/\varepsilon) \Big). \label{eq:bound_m}
\end{equation}

The rest of the proof follows closely the proof of Theorem 3.1 in
\cite{ICML2004_FeatureSelection_Ng}. We outline the key steps for the
sake of completeness. Suppose $m$ satisfies \eqref{eq:bound_m}; then,
with probability at least $1-\delta$, for all $f\in \scrF$
\[ 
\bigg| \frac{1}{m_1}
\sum_{i=1}^{m_1} f(\bphi(x^{(i)}),y^{(i)}) - \mbb{E}_{\bz\sim D} [f(\bz)]
\bigg| \leq \varepsilon.
\]
Thus, using our definition of $\scrF$ in (\ref{classF}), for any $\T$
with $\| \T \|_1 \leq \hat{B}$ and with probability at least $1-\delta$,
we have
\begin{multline*}
\bigg |
\frac{1}{m_1} \sum_{i=1}^{m_1} \big(-\log \mu_{\T}(y^{(i)}|x^{(i)})\big) -
E_{(x,y) \sim D} \big[ -\log \mu_{\T}(y|x) \big]
\bigg |%\\ 
  \leq \varepsilon.
\end{multline*}
Therefore, for all $\T$ with $\| \T \|_1 \leq \hat{B}$ and with
probability at least $1-\delta$, it holds 
\begin{equation*}
| \hat{\epsilon}_{\scrS_1} (\T) - \epsilon(\T) | \leq \varepsilon,
\end{equation*}
where $\scrS_1$ is the training set from Algorithm~\ref{alg:train}. 

Essentially, we have shown that for $m$ large enough the empirical
log-loss function $\hat{\epsilon}_{\scrS_1}(\cdot)$ is a good estimate
of the log-loss function $\epsilon(\cdot)$. According to Step~2 of
Algorithm~\ref{alg:train}, $\hat{\T} = \arg\min_{\T: \|\T\|_1 \leq
  \hat{B}} \hat{\epsilon}_{\scrS_1}(\theta)$. By
Lemma~\ref{lem:4}, we have
\begin{align}
\epsilon(\hat{\T}) \leq & \min_{\T: \|\T\|_1 \leq
  \hat{B}} \epsilon(\T) + 2\varepsilon \notag \\
\leq &  \epsilon(\T^*) + 2\varepsilon, 
\label{truelogloss}
\end{align}
where $\T^*$ is the target policy and the last inequality follows simply
from the fact that $\|\T^*\|_1\leq r K\leq \hat{B}$.  

Eq.~\ref{truelogloss} indicates that the training step of
Algorithm~\ref{alg:train}, finds at least one parameter vector
$\hat{\T}$ whose performance is nearly as good as that of $\T^*$. At the
validation step, we select one of the $\T_B$ found during training.  It
can be shown (\cite{anthony2009neural}) that with a validation set of
the same order of magnitude as the training set (and independent of
$n$), we can ensure that with probability at least $1-\delta$, the
selected parameter vector will have performance at most $2 \varepsilon$
worse than that of the best performing vector discovered in the training
step. Hence, with probability at least $1-2\delta$, the output $\hat{\T}$ of
our algorithm satisfies
\begin{equation}
\epsilon(\hat{\T}) \leq \epsilon(\T^*) + 4 \varepsilon. 
\end{equation}
Finally, replacing $\delta$ with $\delta/2$ and $\varepsilon$ with
$\varepsilon/4$ everywhere in the proof, establishes
Theorem~\ref{thm:1}.

\section{Bounds on Regret} \label{sec:regret}

Theorem~\ref{thm:1} provides a sufficient condition on the number of
samples required to learn a policy whose log-loss performance is close
to the target policy. In this section we study the regret of the
estimated policy, defined in Sec.~\ref{sec: Setup} as the difference
between the average reward of the target policy and the estimated
policy. Given that we use a number of samples in the training set
proportional to the expression provided in Theorem~\ref{thm:1}, we
establish explicit bounds on the regret.

We will bound the regret of the estimated policy by separating the
effect of the error in estimating the policy function (which is
characterized by Theorem~\ref{thm:1}) and the effect the estimated
policy function introduces in the stationary distribution of the Markov
chain governing how states are visited. To bound the regret due to the
perturbation of the stationary distribution, we will use results from
the sensitivity analysis of Markov chains. In
Section~\ref{sec:regret-main} we provide some standard definitions for
Markov chains and state our result on the regret, while in
Section~\ref{sec: Analysis} we provide a proof of this result.

\subsection{Main Result} 
\label{sec:regret-main}

We start by defining the fundamental matrix of a Markov chain.
\begin{defi}
  The fundamental matrix of a Markov chain with state transition
  probability matrix $\BS{P}_{\T}$ induced by RSP $\T$ is
  \[\BS{Z}_{\BS{\T}}=(\BS{A}_{\BS{\T}} +
\BS{e}\BS{\pi'_{\BS{\T}}})^{-1},
\]
where $\BS{e}$ denotes the vector of all $1$'s,
$\BS{A}_{\BS{\T}}=\BS{I}-\BS{P}_{\BS{\T}}$ and $\BS{\pi}_{\BS{\T}}$
denotes the stationary distribution associated with $\BS{P}_{\T}$. Also,
the group inverse of $\BS{A}_{\BS{\T}}$ denoted as $\bA^\#_{\BS{\T}}$ is the
unique matrix satisfying
\[\BS{A}_{\BS{\T}}\BS{A}^\#_{\BS{\T}}\BS{A}_{\BS{\T}}=\BS{A}_{\BS{\T}},
\;
\BS{A}^\#_{\BS{\T}}\BS{A}_{\BS{\T}}\BS{A}^\#_{\BS{\T}}=\BS{A}_{\BS{\T}},
\;
\BS{A}_{\BS{\T}}\BS{A}^\#_{\BS{\T}}=\BS{A}^\#_{\BS{\T}}\BS{A}_{\BS{\T}}.
\]
\end{defi}

Most of the properties of a Markov chain can be expressed in terms of
the fundamental matrix $\BS{Z}$. For example, if $\bpi_1$ is the
stationary probability distribution associated with a Markov chain whose
transition probability matrix is $\BS{P}_1$ and if
$\BS{P}_2=\BS{P}_1+\BS{E}$ for some perturbation matrix $\BS{E}$, then
the stationary probability distribution $\bpi_2$ of the Markov chain
with transition probability matrix $\BS{P}_2$ satisfies the relation
\begin{equation}
\label{eqn:FundamentalMatrix}
\BS{\pi}_1^\prime-\BS{\pi}_2^\prime=\BS{\pi}_2^\prime\BS{EZ_1},
\end{equation}
where $\BS{Z}_1$ is the fundamental matrix associated with $\BS{P}_1$.

\begin{defi}
The ergodic coefficient of a matrix $\BS{B}$ with equal row sums is
\begin{equation}
\label{dfn:ErgodicCoeffi}
\tau(\BS{B})=\sup_{\bv^\prime \BS{e}=0 ; \| \bv \|_1 =1 } \|
\bv^\prime\BS{B}\| _1=\frac{1}{2}\max_{i,j}\sum_{s}
|b_{is}-b_{js}|.
\end{equation}
\end{defi}
The ergodic coefficient of a Markov chain indicates sensitivity of its
stationary distribution. For any stochastic matrix $\BS{P}$, $0\leq
\tau(\BS{P})\leq 1$.

We now have all the ingredients to state our main result bounding regret.
\begin{thm}
\label{thm:RegretBound}
Given $\epsilon>0$ and $\delta>0$, suppose 
\[ m = \Omega((\log n) \cdot
\text{poly} (r, K, C, H,\log(1/\delta), 1/\epsilon))
\] 
i.i.d. samples are
used by Algorithm~\ref{alg:train} to produce an estimate $\hat{\T}$ the
unknown target RSP policy parameter $\T^*$. Then with probability at
least $1-\delta$, we have
\[
|\overline{R}(\T^*)-\overline{R}(\hat{\T})| \leq \sqrt{2\epsilon \log 2 }
R_{\max} ( 1+ \kappa).
\]
where 
\[
R_{\max}= \max_{(x,a)\in \mathcal{X}\times \mathcal{A}} |R(x,a)|
\] 
and $\kappa$ is a constant that depends on the RSP $\hat{\BS{\theta}}$
and can be any of the following:
\begin{itemize}
\item $\kappa=\| \BS{Z}_{\hat{\BS{\T}}} \|_1$,
\item $\kappa= \| \BS{A}^\#_{\hat{\BS{\T}}} \|_1$,
\item $\kappa = 1/(1-\tau(\BS{P}_{\hat{\BS{\T}}} ))$,
\item $\kappa= \tau(\BS{Z}_{\hat{\BS{\T}}}
  )=\tau(\BS{A}^\#_{\hat{\BS{\T}}})$.
\end{itemize}
\end{thm}
The constant $\kappa$ is referred to as {\em condition number}. The
regret is thus governed by the condition number of the estimated RSP;
the smaller the condition number of the trained policy, the smaller is
the regret.

\section{Proof of the Main Result}
\label{sec: Analysis}

In this section we analyze the average reward obtained by the MDP when
we apply the estimated RSP $\hat{\T}$, and prove the regret bound in
Theorem \ref{thm:RegretBound}.  First, we bound the regret as the sum of
two parts.
\begin{eqnarray}
\nonumber
\lefteqn {\mathrm{Reg}(\scrS(\T^*))=\overline{R}(\T^*)-\overline{R}(\hat{\T})}\\
\nonumber
&=&\sum_x \sum_a [\eta_{\T^*}(x,a)-\eta_{\hat{\T}}(x,a)] R(x,a)\\
\nonumber
&=&\sum_x \pi_{\T^*}(x)\sum_a \mu_{\T^*}(a|x)R(x,a)\\
\nonumber
&&-\sum_x \pi_{\hat{\T}}(x)\sum_a \mu_{\hat{\T}}(a|x)R(x,a)\\
\nonumber
&=&\sum_x \pi_{\T^*}(x)\sum_a [\mu_{\T^*}(a|x)-\mu_{\hat{\T}}(a|x)]R(x,a)\\
\nonumber
&&-\sum_x [\pi_{\hat{\T}}(x)- \pi_{\T^*}(x)]\sum_a \mu_{\hat{\T}}(a|x)R(x,a)\\
\nonumber
&\leq & \left | \sum_x \pi_{\T^*}(x)\sum_a [\mu_{\T^*}(a|x)-\mu_{\hat{\T}}(a|x)]R(x,a) \right |\\
\label{eqn:AbsoluteSum}
&&+ \left |\sum_x [\pi_{\hat{\T}}(x)- \pi_{\T^*}(x)]\sum_a \mu_{\hat{\T}}(a|x)R(x,a) \right |.
\end{eqnarray}
Note that the first absolute sum above has terms $\sum_a
[\mu_{\T^*}(a|x)-\mu_{\hat{\T}}(a|x)]$ for all $x$ that are related to
the estimation error from fitting the RSP policy $\hat{\T}$ to
$\T^*$. The second part has terms $ \sum_x | \pi_{\hat{\T}}(x)-
\pi_{\T^*}(x)|$ that are related to the perturbation of the stationary
distribution of the Markov chain by applying the fitted RSP $\hat{\T}$
instead of the original $\T^*$. In the following, we bound each term
separately. We begin with the first term.

\begin{eqnarray}
\nonumber
\lefteqn{  \left | \sum_x \pi_{\T^*}(x)\sum_a
[\mu_{\T^*}(a|x)-\mu_{\hat{\T}}(a|x)]R(x,a) \right | }\\
\nonumber
&\leq & \sum_x \pi_{\T^*}(x)\sum_a |\mu_{\T^*}(a|x)-\mu_{\hat{\T}}(a|x)|
\cdot |R(x,a)|\\
\label{eqn: CS1}
&\leq& R_{\max} \sum_x \pi_{\T^*}(x) \|
\BS{\mu}_{\T^*}(\cdot|x)-\BS{\mu}_{\hat{\T}}(\cdot|x) \|_1, 
\end{eqnarray}
where $\BS{\mu}_{\T}(\cdot|x)$ denotes the probability distribution (a
vector) on the action space $\scrA$ induced by the RSP $\T$ at state
$x$.

The bound in (\ref{eqn: CS1}) is related to the difference in the
log-loss between the RSPs $\T^*$ and $\hat{\T}$. To see this, we need
the following result that connects the $\ell_1$ distance between two
distributions with their KL-divergence.  Let $\BS{p}_1$ and $\BS{p}_2$
denote two probability vectors on $\mathcal{A}$. From the variation
distance characterization of $\BS{p}_1$ and $\BS{p}_2$ we have the
following lemma.
\begin{lemma} (\cite[Lemma 11.6.1]{Book_ElementsOfInformation})
\label{lma:DivergenceInequality}
\begin{eqnarray}
\label{eqn:VariationalDistance}
D(\BS{p}_1 \| \BS{p}_2) \geq \frac{1}{2 \log 2}\| \BS{p}_1-\BS{p}_2\|_1^2.
\end{eqnarray}
\end{lemma}
\vspace{.5cm}

Continuing the chain of inequalities from (\ref{eqn: CS1}), we obtain
\begin{eqnarray}
\nonumber
\lefteqn{  \left | \sum_x \pi_{\T^*}(x)\sum_a [\mu_{\T^*}(a|x)-\mu_{\hat{\T}}(a|x)]R(x,a) \right | }\\
\label{eqn:Variational}\nonumber
&\leq & R_{\max} \sum_x \pi_{\T^*}(x) \sqrt{2 \log 2 D \left (\BS{\mu}_{\T^*}(\cdot|x)\| \BS{\mu}_{\hat{\T}}(\cdot|x) \right )} \\
\nonumber
&\leq & \sqrt{2 \log 2} R_{\max} \cdot\\
\label{eqn:Jenson1}\nonumber
&&  \sqrt{\sum_x \pi_{\T^*}(x) D \left (\BS{\mu}_{\T^*}(\cdot|x)\| \BS{\mu}_{\hat{\T}}(\cdot|x) \right )}\\
\label{eqn:Bound1}
&= & \sqrt{2 \log 2}R_{\max} \sqrt{D_{\T^*}(\BS{\mu}_{\T^*}\|
  \BS{\mu}_{\hat{\T}})}. 
\end{eqnarray}
In the first inequality, we applied Lemma \ref{lma:DivergenceInequality}
by setting $\BS{p}_1=\BS{\mu}_{\T^*}(\cdot|x)$ and
$\BS{p}_2=\BS{\mu}_{\hat{\T}}(\cdot|x)$ for each $x$. In second
inequality, we applied Jensen's inequality. We can now use Theorem
\ref{thm:1} to bound (\ref{eqn:Bound1}).

We next bound the second term in (\ref{eqn:AbsoluteSum}) using
techniques from perturbation analysis of eigenvalues of a matrix. We
have 
\begin{eqnarray}
\nonumber
\lefteqn{\left | \sum_x (\pi_{\hat{\T}}(x)-\pi_{\T^*}(x))\sum_a \mu_{\hat{\T}}(a|x)R(x,a)\right |}\\
\nonumber
&\leq & \sum_x |\pi_{\hat{\T}}(x)-\pi_{\T^*}(x)|\sum_a  |\mu_{\hat{\T}}(a|x) R(x,a)| \\
\label{eqn:CS2}
&\leq & R_{\max}\sum_x |\pi_{\hat{\T}}(x)-\pi_{\T^*}(x) | \sum_{a}
\mu_{\hat{\T}}(a|x) \notag \\
&= & R_{\max}\sum_x | \pi_{\hat{\T}}(x)-\pi_{\T^*}(x)|,
\label{eqn:ProbabilityNorm}
\end{eqnarray}
where we used the definition of $R_{\max}$ in the 2nd inequality and the
last equality follows by noting that $ \sum_a \mu_{\hat{\T}}(a|x) = 1$
for all $x$.  To further bound the difference in stationary
distributions in (\ref{eqn:ProbabilityNorm}), we use a relation between
a measure of perturbation of the stationary distribution and the
condition number of Markov chains. We recall the following result that
is useful to bound the difference between the stationary distribution
induced by the optimal RSP and that induced by the estimated RSP in
terms of the condition number of the Markov chains.

\begin{lemma}[\cite{APT_PerturbationOfStationary_Seneta,
Elsevier01_ComparisonOfPerturbation_ChoMeyer}]
\label{lemma:ErgodicCOeff}
Let $\BS{\pi}_{\BS{\T}_1}$ and $\BS{\pi}_{\BS{\T}_2}$ be the unique
stationary distributions of the stochastic matrices $\BS{P}_{\BS{\T}_1}$
and $\BS{P}_{\BS{\T}_2}$, respectively. Let
$\BS{E}=\BS{P}_{\BS{\T}_1}-\BS{P}_{\BS{\T}_2}$. Then,
\begin{equation}
\label{eqn:PerturbationDistribution}
\| \BS{\pi}_{\BS{\T}_1} -\BS{\pi}_{\BS{\T}_2} \|_1 \leq \kappa \|
\BS{\pi}_{\T_1}'\BS{E}\|_1,
\end{equation}
where $\kappa$ is a constant that can take the following values
\begin{itemize}
	\item $\kappa=\| \BS{Z}_{\BS{\T}_2} \|_1$,
	\item $\kappa= \| \BS{A}^\#_{\BS{\T}_2}\|_1$,
	\item $\kappa = 1/(1-\tau(\BS{P}_{\BS{\T}_2}))$,
	\item $\kappa= \tau(\BS{Z}_{\BS{\T}_2})=\tau(\BS{A}^\#_{\BS{\T}_2})$.
\end{itemize}
\end{lemma}
\pf
  The proof follows by setting $\BS{\pi}_1=\BS{\pi}_{\T_2}$ and
  $\BS{\pi}_2=\BS{\pi}_{\T_1}$ in (\ref{eqn:FundamentalMatrix}) and
  using the relation
\[
\| \BS{\pi}_{\BS{\T}_1} -\BS{\pi}_{\BS{\T}_2} \|_1 = \|
\BS{\pi}_{\T_1}' \BS{EZ}_{\T_2} \| _1\leq \|
\BS{\pi}_{\T_1}' \bE \|_1 \| \BS{Z}_{\T_2} \|_1.
\] 
The other relations follow similarly from Sec. 3 of
\cite{Elsevier01_ComparisonOfPerturbation_ChoMeyer}.
\qed

Continuing the chain of inequalities in (\ref{eqn:ProbabilityNorm}) and
applying Lemma~\ref{lemma:ErgodicCOeff} by setting
$\pi_{\BS{\T}_1}=\pi_{\T^*}$ and $\pi_{\BS{\T}_2}=\pi_{\hat{\T}}$, we obtain
\begin{multline}
\label{eqn:ErgodicCoefficient}
\left| \sum_x (\pi_{\hat{\T}}(x)-\pi_{\T^*}(x))\sum_a
\mu_{\hat{\T}}(a|x)R(x,a)\right | %\\
\leq  R_{\max}
\kappa \|\bpi_{\T^*}' (\BS{P}_{\T^*}-\BS{P}_{\hat{\T}})\|_1,
\end{multline}
where, with some overloading of notation, $\kappa$ is now as specified
in the expressions provided in the statement of
Theorem~\ref{thm:RegretBound}.

The $i$th component of the vector $\BS{\pi}_{\T^*}'
(\BS{P}_{\T^*}-\BS{P}_{\hat{\T}})$ is given by
\begin{eqnarray*}
\lefteqn{[\BS{\pi}_{\T^*}' (\BS{P}_{\T^*}-\BS{P}_{\hat{\T}})]_i}\\
&=&  \sum_x \pi_{\T^*}(x)[\BS{P}_{\T^*}(i|x)-\BS{P}_{\hat{\T}}(i|x)] \\
&=& \sum_x
\pi_{\T^*}(x)\sum_a[\BS{P}(i|x,a)(\mu_{\T^*}(a|x)-\mu_{\hat{\T}}(a|x))],
\end{eqnarray*}
where, in the last equality, we applied the definition of the transition
probability $\BS{P}_{\T}$ associated with RSP $\T$. 

It follows
\begin{align}
\nonumber
& \| \BS{\pi}^\prime_{\T^*}(\BS{P}_{\T^*}-\BS{P}_{\hat{\T}})\|_1\\
\nonumber
=& \sum_y \left |\sum_x \pi_{\T^*}(x)\sum_a[\BS{P}(y|x,a)(\mu_{\T^*}(a|x)-\mu_{\hat{\T}}(a|x))]\right |\\
\nonumber
\leq & \sum_x \pi_{\T^*}(x)\sum_y\sum_a\left |[\BS{P}(y|x,a)(\mu_{\T^*}(a|x)-\mu_{\hat{\T}}(a|x))]\right |\\
\label{eqn:PerturbationBound}
\leq & \sum_x \pi_{\T^*}(x)\sum_a \left |
\mu_{\T^*}(a|x)-\mu_{\hat{\T}}(a|x)\right |,
\end{align}
where the last inequality follows by noting that
$\sum_{y}\BS{P}(y|x,a)=1$ for all $(x,a)$. 

Now, using (\ref{eqn:PerturbationBound}), (\ref{eqn:ProbabilityNorm}),
similar steps as in (\ref{eqn:Bound1}), and Theorem~\ref{thm:1}, we can
bound the second term in (\ref{eqn:AbsoluteSum}) as
\begin{multline}
\label{eqn:Bound2}
\left | \sum_x (\pi_{\hat{\T}}(x)-\pi_{\T^*}(x))\sum_a
  \mu_{\hat{\T}}(a|x)R(x,a)\right | %\\
\leq  \sqrt{2\epsilon \log 2 }\ \kappa\ R_{\max}.
\end{multline}
Finally, combining (\ref{eqn:Bound1}) and (\ref{eqn:Bound2}) and
applying Theorem~\ref{thm:1}, the result in
Theorem~\ref{thm:RegretBound} follows.

\section{A robot navigation example}\label{sec: Simulation}

In this section, we discuss the experimental setup to simulate an MDP
and validate the effectiveness of our proposed learning algorithm. We
consider the problem of learning the policy used by an agent (a robot)
as it moves on a 2-dimensional grid.  After simulating the movement of
the robot and recording its actions in various states, we use these
states-action samples to learn the policy of the robot and evaluate its
performance.

\subsection{Environment and Agent Settings}

\begin{figure}[!ht]
\centering
\includegraphics[width=0.7\textwidth]{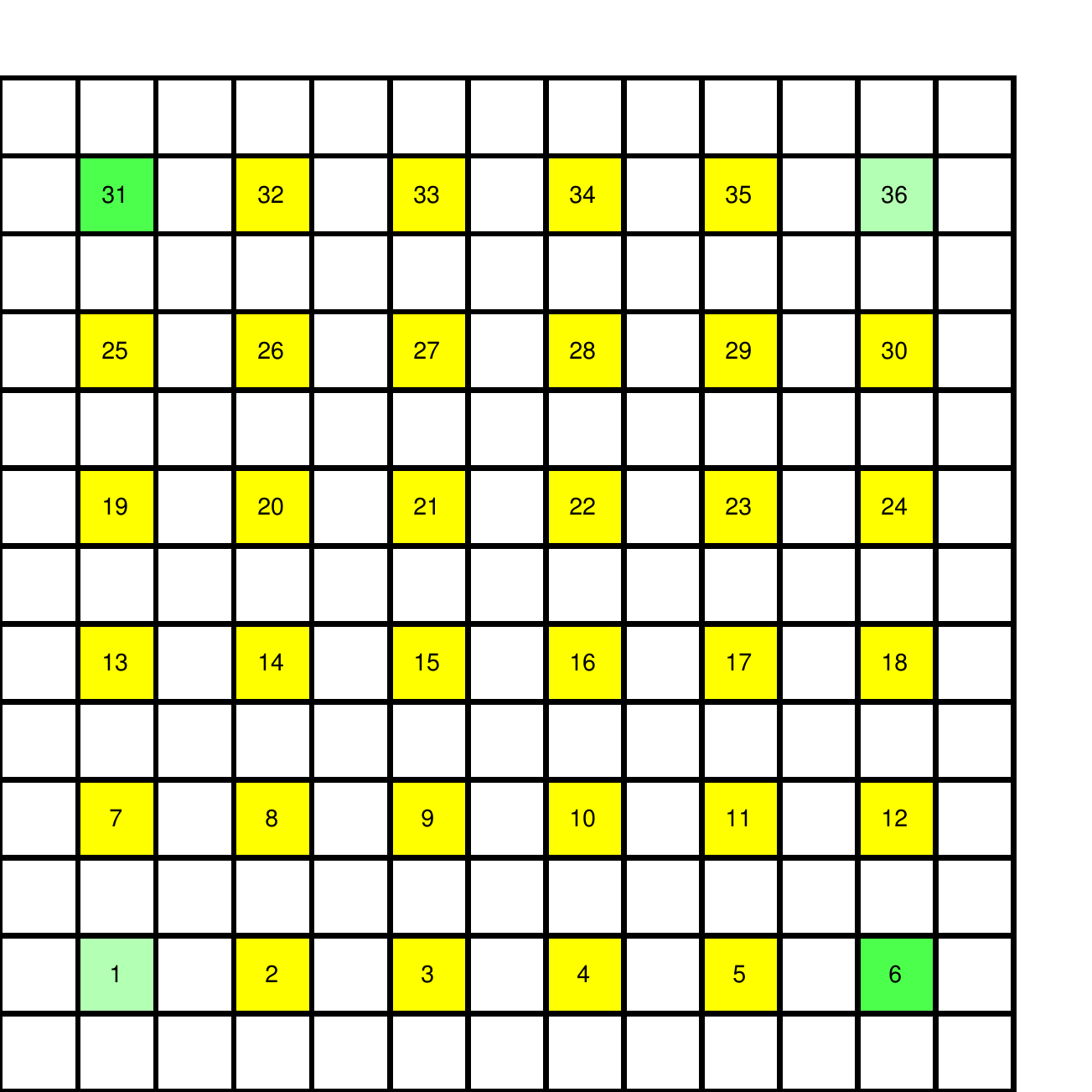}
\caption{The environment of the MDP is a $13 \times 13$ grid. Colored
  (or shaded) grid squares correspond to waypoints (identified by the
  southwest vertex of the square) for defining features that are used by
  the estimated policy. The four (extreme) squares labeled $1$, $6$,
  $31$, $36$ correspond to waypoints on the grid with associated
  reward. The one step reward function of the MDP is a weighted sum of
  Gaussian functions specified by these reward points on the grid.}
\label{fig:environment}
\end{figure}

Consider an agent moving in a $13\times 13$ grid, shown in
Fig. \ref{fig:environment}. The agent's position is specified by a
two-tuple {\em state} $\bx=(x_1, x_2)$, representing its coordinates in
the grid.  We assume $(0, 0)$ is at the southwest corner of the grid and
we make the convention that coordinates $\bx=(x_1, x_2)$ on the grid
identify the square defined by the four points $(x_1, x_2)$, $(x_1+1,
x_2)$, $(x_1, x_2+1)$, and $(x_1+1, x_2+1)$. For example, when we say
that the agent is at $\bx=(x_1, x_2)$ we mean that the agent can be
anywhere in the above square.

At each time instance, the agent can take $4$ actions: \texttt{North},
\texttt{East}, \texttt{West}, and \texttt{South}. Without loss of
generality, we assume that the agent can only move to a neighboring grid
point. The destination of the agent is based on its current position and
the action. We assume that the agent movements are subject to
uncertainty, which can cause the agent's intended next position to shift
to a point adjacent to that position. For example, the agent in state
$(x_1, x_2)$ taking action \texttt{North} will enter state $(x_1,
x_2+1)$ with probability $0.8$ (the intended next position), but can
enter states $(x_1-1, x_2)$ or $(x_1+1, x_2)$ with probability $0.1$,
respectively. At the boundary points of the grid, the agent bounces
against the ``wall'' in the opposite direction with its position
unchanged.

\begin{figure}[!ht]
\centering
\includegraphics[width=0.7\textwidth]{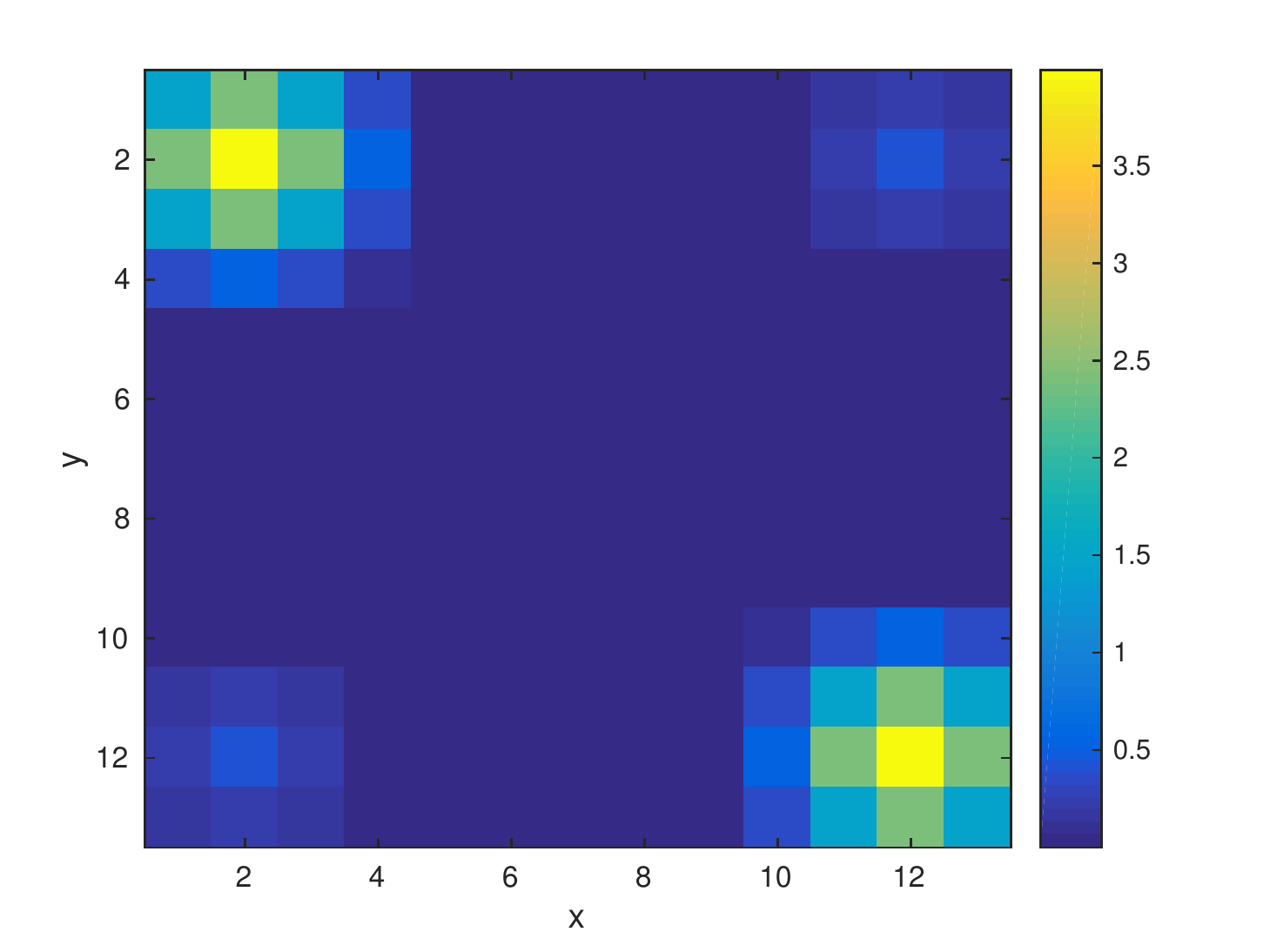}
\caption{The one step reward function induced by the reward waypoints
  shown in Fig. \ref{fig:environment}.}
\label{fig:rewards}
\end{figure}

The environment contains points with associated rewards. Specifically,
as showed in Fig. \ref{fig:environment}, squares (waypoints) labeled $1$
and $36$ have associated reward equal to $1$ and squares (waypoints)
labeled $6$ and $31$ have associated reward equal to $20$. The reward
$r_\bx$ of a waypoint $\bx$ ``spreads'' on the grid according to a
Gaussian function. Specifically, the immediate reward at point $\by$ 
due to the reward waypoint $\bx$ is given by
\[ 
r_\bx \frac{1}{\sqrt{2\pi}} e^{ - \frac{\|\bx-\by\|_2^2}{2}}
\]
and is constant for all actions. Summing over all reward waypoints
$\bx$, the immediate reward at some point $\by$ is given by
\begin{equation}
\label{f-reward}
f(\by) = \sum_\bx r_\bx \frac{1}{\sqrt{2\pi}} e^{ - \frac{\|\bx-\by\|_2^2}{2}},
\end{equation}
where we assume that a point $\bx$ with no associated reward satisfies
$r_\bx=0$. The one step reward function induced by the four reward
waypoints of Fig.~\ref{fig:environment} is shown in
Fig. \ref{fig:rewards}.

The agent (robot) entering some state (position) $\by$ collects reward
equal to the immediate reward $f(\by)$.  The objective of the agent is
to navigate on the grid so as to maximize the long-term average reward
collected.

The key step of efficient learning is to define appropriate features we
will use to represent the agent's policy. The way of selecting features
in this paper is inspired by spline interpolation
\cite{de1978splines}. We define $36$ waypoints on the grid, shown in
Fig.~\ref{fig:environment}. For a state-action pair $(\bx,a)$, we set
$\by$ to be the intended next state and define the features
\begin{equation*}
\phi_i(\bx,a) = f_i(\by),\quad i=1,\ldots,36,
\end{equation*}
where $f_i(\cdot)$ is as defined as
\[
f_i(\by) = r_{\bx_i} \frac{1}{\sqrt{2\pi}} e^{ -
  \frac{\|\bx_i-\by\|_2^2}{2}},
\]
and $\bx_i$ is the location of the $i$th waypoint. 

\subsection{Simulation and Learning Performance}

For the MDP we have introduced, we use value iteration \cite{bert-dp} to
find the best policy. We generate independent state-action samples
according to this policy. Then, we estimate the policy using the above
samples according to the logistic regression algorithm discussed in
Section~\ref{sec: Regression}. We largely follow the style of
\cite{ICML2004_FeatureSelection_Ng} in presenting our numerical
results. We compare average rewards from $3$ different policies with
respect to the number of features and the number of samples as follows:
\begin{enumerate}
\item {\em Target policy}: We choose the policy obtained by the value
  iteration algorithm \cite{bert-dp} as the target policy. It is used to
  generate samples for learning.
	
\item {\em $\ell_1$-regularized policy}: The RSP trained using the
  Algorithm in Section~\ref{sec: Regression}.
	
\item {\em Unregularized policy}: The RSP trained using logistic
  regression in Section~\ref{sec: Regression}, but without the $\ell_1$
  constraint on the parameter vector $\T$ (cf. problem
  (\ref{eqn:logistic})).
	
\item {\em Greedy policy}: The agent takes the action with the largest
  expected next step reward, i.e., the local reward feature is the only
  consideration. This policy is used as a baseline.
\end{enumerate}

\begin{figure}[!ht]
\centering
\includegraphics[width=0.7\textwidth]{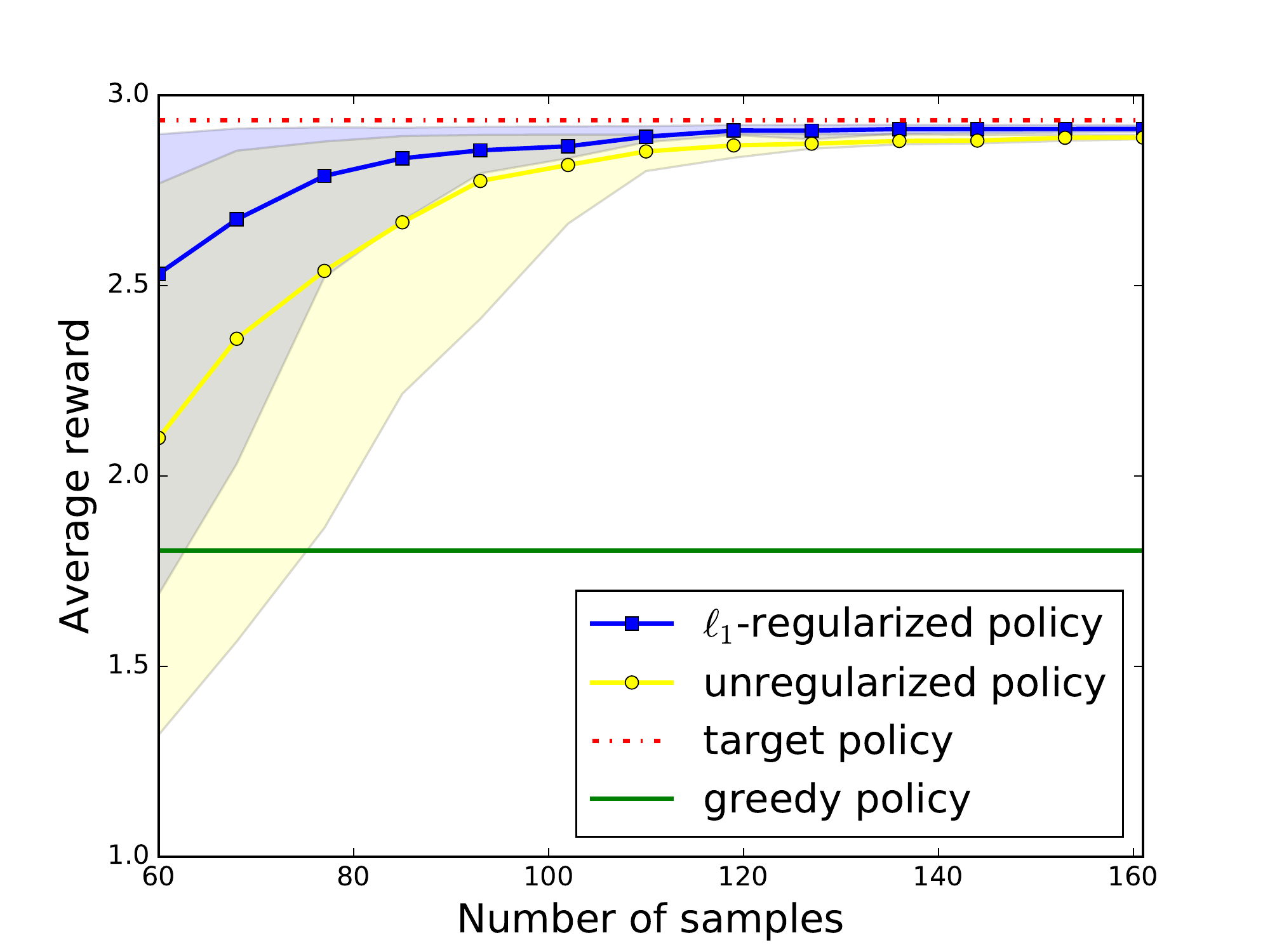}
\caption{Average rewards (over 100 runs) of different policies as a function of the number of samples.}.\label{fig:samples}
\end{figure}

We randomly sample the state-action pairs generated by the optimal
policy and use these samples to form a training set to be used in
learning an RSP (as in Section~\ref{sec: Regression}). This process is
repeated $100$ times. The average rewards of all policies we considered
as a function of the number of samples used for training are displayed
in Fig.~\ref{fig:samples}. As shown in Fig.~\ref{fig:samples}, the
$\ell_1$-regularized policy performs closely to the target policy, even
with a very small number of training samples. Not surprisingly, both the
$\ell_1$-regularized and the unregularized policy perform almost as well
as the optimal policy when the number of training samples grows large.
For most runs, the policies we learn (both regularized and
unregularized) perform better than the greedy policy. 

It is interesting that when the number of samples is small, the
regularized policy performs better. This is because the unregularized
policy tends to overfit to the training samples while the regularized
policy seeks to learn fewer parameters and thus ends up learning them
better, in a way that generalizes well out-of-sample. In sum, this
simulation example supports the effectiveness of our algorithm.

\section{Conclusion}\label{sec:Conclusion}

We considered the problem of learning a policy in a Markov decision
process using the state-action samples associated with the policy. We
focused on a Boltzmann-type policy that is characterized by feature
vectors associated with each state-action pair and a parameter that is
sparse. 

To learn the policy, we used $\ell_1$-regularized logistic regression
and showed that a good generalization error bound also guarantees a good
bound on the regret, defined as the difference between the average
reward of the estimated policy and the target policy. Our results
suggest that one can estimate an effective policy using a training set
of size proportional to the logarithm of the number of features used to
represent the policy.

%\ifCLASSOPTIONcaptionsoff
%\newpage
%\fi


%\bibliographystyle{IEEEtran}
%\bibliography{/home/yannisp/Private/bib/abbrev,/home/yannisp/Private/bib/IEEEabrv,/home/yannisp/Private/bib/communications,/home/yannisp/Private/bib/my,/home/yannisp/Private/bib/optimization,/home/yannisp/Private/bib/stochastics,/home/yannisp/Private/bib/various,/home/yannisp/Private/bib/scheduling,Learning}

%\bibliography{C:/Users/yannisp/Documents/Private/bib/IEEEabrv,C:/Users/yannisp/Documents/Private/bib/abbrev,C:/Users/yannisp/Documents/Private/bib/communications,C:/Users/yannisp/Documents/Private/bib/my,C:/Users/yannisp/Documents/Private/bib/optimization,ano,cyber,web}


\end{document}